\documentstyle[12pt]{article}
\begin{document}

\begin{center}
{\bf ANALOGS OF {\boldmath$q$}-SERRE RELATIONS IN THE YANG-BAXTER ALGEBRAS}
\end{center}

\vspace*{-.1cm}

\begin{center}
{\sc M. L\"udde} \\[.1cm]
{\em SAP AG, Basis 02, Postfach 1461,
D-69185 Walldorf/Baden, Germany \\ Mirko.Luedde@SAP-AG.De}
\end{center}
\begin{center}
{\sc A.\,A. Vladimirov} \\[.1cm]
{\em BLTP, JINR, Dubna, Moscow region 141980, Russia \\
alvladim@thsun1.jinr.ru}
\end{center}

\vspace{.7cm}

{\small 
Yang-Baxter bialgebras,
as previously introduced by the authors,
are shown to arise
from a double crossproduct construction applied to the bialgebra
$$
R_{12}T_1T_2=T_2T_1R_{12}\,, \ \ E_1T_2=T_2E_1R_{12}\,, \ \
\Delta(T)=T\hat{\otimes} T\,, \ \ \Delta(E)=E\hat{\otimes}
T+1\hat{\otimes} E
$$
and its skew dual, with $R$ being a numerical matrix solution of the
Yang-Baxter equation.
It is further shown that a set of relations
generalizing $q$-Serre ones in the Drinfeld-Jimbo algebras $U_q({\bf
g})$
can be naturally imposed on Yang-Baxter algebras from the requirement
of non-degeneracy of the pairing.
}

\vspace{.6cm}

{\bf 1}.
Yang-Baxter algebras (YBA),
introduced in \cite{Lud,Vl1,Lud2},
generalize the wide-known FRT construction~\cite{FRT}
in the following sense: to any numerical matrix solution $R$ of the
Yang-Baxter equation there is associated a bialgebra containing the FRT
one
as a sub-bialgebra.  Generally, this construction may provide examples
of
(new) bialgebras and Hopf algebras~\cite{Vl2}. In several
aspects, there is some similarity of YBA with so-called inhomogeneous
quantum groups~\cite{SWW,Ma1}, and also with Majid's scheme of double
bosonization~\cite{Ma2}. However, in YBA no extra (dilation)
generators appear, whereas additional (analogous to $q$-Serre) relations
are not necessarily quadratic. The main goal of the present note is to
refine the concept and definition of these generalized Serre relations,
first introduced in~\cite{Lud} and further studied in~\cite{Vl2}.
Here we obtain Serre-like relators as the elements
in the kernel of some bilinear form, and the main result of the present
paper is a condition (\ref{1.19}) for such relators.

\vspace{.2cm}

{\bf 2}.
We recall the definition of the Yang-Baxter bialgebra.
Let an invertible matrix $R$ obey the Yang-Baxter equations
$R_{12}R_{13}R_{23}=R_{23}R_{13}R_{12}$.
Consider two bialgebras
$ Y_+ $ and $ Y_- $
with generators
$ \{ u_{j}^{i}, F^i \} , $
$ \{ t_{j}^{i}, E_i \} , $
respectively,
which form matrices
$T,U$, a row $E$ and a column $F$.
We impose the following multiplication,
\begin{equation}
\label{1} \begin{array}{ll}
R_{12}\,T_1\,T_2=T_2\,T_1\,R_{12}\,,
&
E_1\,T_2=T_2\,E_1\,R_{12}\,,
\\
R_{12}\,U_1\,U_2=U_2\,U_1\,R_{12}\,,
&
F_2\,U_1=R_{12}\,U_1\,F_2\,,
\end{array}
\end{equation}
comultiplication and counit, respectively:
\begin{equation}
\label{2}  \begin{array}{llll}
\Delta(T)=T\hat{\otimes} T\,, & \Delta(E)=E\hat{\otimes}
T+1\hat{\otimes} E\,,
& \varepsilon (T)={\bf 1}\,, & \varepsilon (E)=0\,, \\
\Delta(U)=U\hat{\otimes} U\,, & \Delta(F)=F\hat{\otimes}
1+U\hat{\otimes} F\,,
& \varepsilon (U)={\bf 1}\,, & \varepsilon (F)=0.
\end{array}
\end{equation}
The symbol $\hat{\otimes}$ denotes the tensor product with
implied matrix multiplication.
There is a bilinear pairing
$ Y_+ \otimes Y_- \rightarrow {\bf C} $
satisfying
\begin{equation}
\label{3}
<U\,,T>=R\,,\ \ \ <1\,,T>=<U\,,1>=<F\,,E>={\bf 1}
\end{equation}
such that these are the only nonzero pairings between generators
and such that $ Y_+ $ and $ Y_- $
are `skew paired'~\cite{Ma4}:
\begin{equation}
\label{4}
<xy\,,a>=<x\otimes y\,,\Delta(a)>\,, \ \
<\Delta(x)\,,a\otimes b>=<x\,,ba>\,,
\end{equation}
\begin{equation}
\label{5}
\varepsilon (a)=<1,a>\,, \ \ \ \varepsilon (x)=<x,1>\,.
\end{equation}

\vspace{.1cm}

{\bf 3}.
The bialgebras defined above form a matched pair of bialgebras,
when a left action of $ Y_+ $ on $ Y_- $
and a right action of $ Y_- $ on $ Y_+ $ are suitably defined.
So they can be used to build a `double crossproduct'~\cite{Ma4,Ma5}
bialgebra
$ Y. $
This is because there is another pairing
\begin{equation}
\label{6}
<U\,,T>^- = R^{-1}\,,\ \ \
<1\,,T>^- = <U\,,1>^- = {\bf 1}\,, \ \ <F\,,E>^- = {\bf -1},
\end{equation}
yielding a `convolution inverse' of the pairing (\ref{3}).
This means,
\begin{equation}
\label{7}
<x_{(1)}\,,a_{(1)}>^- <x_{(2)}\,,a_{(2)}>
=
<x_{(1)}\,,a_{(1)}><x_{(2)}\,,a_{(2)}>^-
=
\varepsilon (x) \varepsilon (a)\,.
\end{equation}
In the Hopf algebra case, $<x,a>^-$ is simply $<S(x),a>$.
For details, see~\cite{Ma4}.
The double crossproduct algebra $ Y $ implies additional
`cross-multiplication' relations
\begin{equation}
\label{ML1} \begin{array}{ll}
R_{12}\,U_1\,T_2=T_2\,U_1\,R_{12}\,,
&
T_2\,F_1=R_{12}\,F_1\,T_2\,,
\\
E\,F-F\,E=T-U\,,
&
U_1\,E_2=E_2\,U_1\,R_{12}\,,
\end{array}
\end{equation}
and $ Y_+ $ and $ Y_- $ are imbedded as sub-bialgebras.
In~\cite{Vl1} it is shown how an antipode can be introduced into this
bialgebra, in \cite{Lud2} a vector space representation is constructed.

\vspace{.1cm}

{\bf 4}.
Commutation relations between the $E_i$\,-generators
themselves,
as well as between the $F^i$, are still missing in
(\ref{1}),(\ref{ML1}).
In~\cite{Lud,Vl2}
two slightly different
(but possibly equivalent)
recipes have been
proposed how to add such relations without destroying the
bialgebra (respectively the Hopf algebra) structure.
The first
of them~\cite{Lud,Lud2} was motivated by geometrical ideas and by the
desire to obtain, via appropriate choices of the matrix $ R, $
the $ U_q({\bf g}) $ algebras.
The second~\cite{Vl2} approach
was in fact an explicit solution of a bialgebra
condition.
However, in the present work the authors choose to
propose one more approach to the problem.
Now we exploit an idea which
is well known in the Lie algebra theory,
see e.g.~\cite{Kac},
and has been used in \cite{Lus} to recover $q$-Serre identities in
$U_q({\bf g})$: Serre relators nullify (lie in the kernel of) a suitable
bilinear form.

\vspace{.2cm}

{\bf 5}.
Let us first recall a general fact about null ideals~\cite{Ma5}.
Let the bialgebras $A$ and $B$ be (skew) paired and let
$I := \{a\in A\,:\,<B\,,a>=0\}$.
Then $I$ is an ideal in $A$ (`null ideal'), and can be shown to be
also a bi-ideal.
Indeed,
let $\{e_i\,,e_\alpha \}$ and $\{e_i\}$ be bases of
$A$ and $I$, respectively. Then, due to (\ref{4}),
\begin{equation}
\label{9}
<B\,,e_i>=0\,, \ \ \
<B\otimes B\,,\Delta(e_i)>=0\,.
\end{equation}
One can write down
\begin{equation}
\label{11}
\Delta(e_i)=f_{i}^{jk}(e_j\otimes e_k)+
f_{i}^{\alpha k}(e_\alpha \otimes e_k)+
f_{i}^{j \beta}(e_j\otimes e_\beta)+
f_{i}^{\alpha \beta}(e_\alpha \otimes e_\beta )\,.
\end{equation}
To prove that $f_{i}^{\alpha\beta}=0$
we choose a set
$\{e^\alpha\}\in B$ such that
$<e^\alpha\,,e_\beta>=\delta_\beta^\alpha$.
Then
\begin{equation}
\label{13}
0=<e^\alpha\otimes e^\beta\,,\Delta(e_i)>=f_{i}^{\mu\nu}
<e^\alpha\otimes e^\beta\,,e_\mu\otimes e_\nu>=f_{i}^{\alpha\beta}\,.
\end{equation}
Consequently, $I$ is a bi-ideal in the bialgebra $A$:
\begin{equation}
\label{14}
\Delta(I)\subset I\otimes A+A\otimes I\,.
\end{equation}
Analogously, $I$ is shown to be a Hopf ideal in the case of
paired Hopf algebras.

\vspace{.2cm}

{\bf 6}.
The $\{U\}$- and $\{T\}$-subalgebras of $ Y_+ $ and $ Y_-, $
respectively,
are bialgebras that are (skew) paired by (\ref{3}).
For general $R$, this pairing may be degenerate, so that `pure-$U$' and
`pure-$T$' null bi-ideals $I_+$ and $I_-$ may exist.
Assume that we know both of these bi-ideals
for a given $R$.
Now let us extend these ideals to ideals
$ I_+' $ and $ I_-' $
in the full bialgebra $Y$ by
multiplying them from left and right, say
$$
I \rightarrow I' := \sum Y \cdot I \cdot Y\,.
$$
These extended ideals are also bi-ideals
due to (\ref{14}) and homomorphisity of the co-product:
$$
\Delta(Y I Y)
= \Delta(Y) \Delta(I) \Delta(Y)
\subset (Y \otimes Y) (I \otimes Y + Y \otimes I) (Y \otimes Y)
\subset I' \otimes Y + Y \otimes I'\,.
$$
Now we divide the algebra $ Y $ by the bi-ideal $ I_+' + I_-' $
and further work with the quotient bialgebra (and the two subalgebras),
keeping the notations as before.
For `regular' $R$~\cite{Ma5},
the double crossproduct structure is respected by the quotient.

\vspace{.2cm}

{\bf 7}.
However, there still might be left some degeneracy
of the pairing between the (quotient) $Y_+$ and $Y_-$ bialgebras.
For example, elements of the form
\begin{equation}
\label{1.9}
E_{a_1}\ldots E_{a_N}\,\omega ^{a_1\ldots a_N}
\end{equation}
with $\omega $ being some numerical
coefficients could also nullify the bilinear form (\ref{3}).
Those elements (\ref{1.9})
(and their analogs formed by $F$-generators) which possess this property
are precisely the analogs of $q$-Serre relators we are interested in.
Now our goal is to find all such elements and then set them to zero,
thus introducing Serre-like relations on $E$- and $F$-generators.

\vspace{.2cm}

{\bf 8}.
Let us find out when an element of type (\ref{1.9}) nullifies the
pairing with all elements of the $Y_+$-bialgebra.
Consider first the case $N=2$, and begin with the following computation:
$$
<F^mF^n\,,E_iE_j\omega ^{ij}>=<\Delta(F^m)\Delta(F^n)\,,
E_j\otimes E_i>\omega ^{ij}
$$
$$
=<(F^m\otimes 1+u_p^m\otimes F^p)(F^n\otimes 1+u_q^n\otimes F^q)\,,
E_j\otimes E_i>\omega ^{ij}
$$
$$
=<F^mu^n_q\otimes F^q+u^m_pF^n\otimes F^p\,,
E_j\otimes E_i>\omega ^{ij}=
<F^mu^n_i+u^m_iF^n\,,E_j>\omega ^{ij}
$$
$$
=(<F^m\otimes u^n_i\,,E_k\otimes t^k_j>+<u^m_i\otimes F^n\,,1\otimes
E_j>)
\omega ^{ij}=(<u^n_i\,,t^m_j>+\delta^m_i\delta^n_j)\omega ^{ij}
$$
\begin{equation}
\label{1.10}
=(R_{ij}^{nm}+{\bf 1}_{ij}^{mn})\omega ^{ij}=
(B+{\bf 1})_{ij}^{mn}\omega^{ij}=([2!]_B\,\omega)^{mn}\,.
\end{equation}
We used here the standard definition of the `braid matrix' $B$:
\begin{equation}
\label{1.7}
B_{mn}^{ij}:=R_{mn}^{ji}\,, \ \ B_m:=B_{m,m+1}\,,
\ \ B_1B_2B_1=B_2B_1B_2\,.
\end{equation}
Definition and properties of `braided factorial' $[N!]_B$ can be found
in~\cite{Lud,Lud2,Ma6}.

Now, consider a general element $F^mF^n\Phi(u)$ of the bialgebra $Y_+$,
which is relevant in the $N=2$ case
(here $\Phi(u)$ may, of course, carry its own indices):
$$
<F^mF^n\Phi(u)\,,E_iE_j>\omega ^{ij}=
<F^mF^n\otimes\Phi(u)\,,\Delta(E_i)\Delta(E_j)>\omega ^{ij}
$$
$$
=<F^mF^n\,,E_pE_q><\Phi(u)\,,t^p_it^q_j>\omega^{ij}
=<\Phi(u)\,,({\bf 1}+B)_{pq}^{mn}\,t^p_it^q_j>\omega^{ij}
$$
\begin{equation}
\label{1.11}
=<\Phi(u)\,,t^m_pt^n_q>({\bf 1}+B)_{ij}^{pq}\omega^{ij}
=<\Phi(u)\,,t^m_pt^n_q>([2!]_B\,\omega)^{pq}\,,
\end{equation}
where the relation
\begin{equation}
\label{1.12}
B_{12}T_1T_2=T_1T_2B_{12}
\end{equation}
has been used. Thus, we obtain
\begin{equation}
\label{1.13}
[2!]_B\,\omega=0
\end{equation}
as a necessary and sufficient condition for $ E_iE_j\,\omega^{ij} $
to lie in the null ideal.

In the general case of arbitrary $N$ in (\ref{1.9}), we find:
$$
<F^{q_1}\ldots F^{q_N}\,,E_{a_1}\ldots E_{a_N}>
=<\Delta(F^{q_1})\ldots\Delta(F^{q_N})\,,
E_{a_N}\otimes E_{a_1}\cdot\ldots\cdot E_{a_{N-1}}>
$$
$$
=\sum_{k=1}^{N}<\!u_{p_1}^{q_1}\ldots u_{p_{k-1}}^{q_{k-1}}
F^{q_k}u_{p_{k+1}}^{q_{k+1}}\ldots u_{p_N}^{q_N},E_{a_N}\!>
<\!\ldots F^{p_{k-1}}F^{p_{k+1}}\ldots,
E_{a_1}\ldots E_{a_{N-1}}\!>
$$
\begin{equation}
\label{1.16}
=\sum_{k=1}^{N}B_{p_kb_{k+1}}^{q_kq_{k+1}}
B_{p_{k+1}b_{k+2}}^{b_{k+1}q_{k+2}}\ldots B_{p_{N-1}a_N}^{b_{N-1}q_N}
<\!F^{q_1}\ldots F^{q_{k-1}}F^{p_k}\ldots F^{p_{N-1}},
E_{a_1}\ldots E_{a_{N-1}}\!>\,.
\end{equation}
By induction, using recursion formulas for braided
factorials~\cite{Lud,Ma6}, we come to
\begin{equation}
\label{1.17}
<F^{q_1}\ldots F^{q_N}\,,E_{a_1}\ldots E_{a_N}>
=([N!]_B)_{a_1\ldots a_N}^{q_1\ldots q_N}
\end{equation}
(compare similar formulas in~\cite{Ma2}).
So, the null ideal condition
\begin{equation}
\label{1.19}
[N!]_B\omega =0
\end{equation}
is proved, because a proper generalization $F\ldots F\Phi(u)$
is easily handled in complete analogy with (\ref{1.11}).

Quite analogously, it can be shown that elements of the form
\begin{equation}
\label{1.9a}
\eta_{a_1\ldots a_N}\,F^{a_1}\ldots F^{a_N}
\end{equation}
lie in the kernel of (\ref{3}) if and only if
\begin{equation}
\label{1.19a}
\eta\,[N!]_B =0\,.
\end{equation}

{\bf 9.}
We now consider the most
general null element containing $E$ generators.
Explicitly using the relation $ET=TER$,
we can cast it to the form
\begin{equation}
\label{2.1}
\omega^{[\alpha]}(T)E_{[\alpha]}
:=
\omega^{\alpha_1\ldots\alpha_N}(T)E_{\alpha_1}\ldots E_{\alpha_N}\,,
\end{equation}
allowing polynomial dependence of the coefficients
$\omega $ on the generators $T$.
Then a condition of being in the kernel
is derived as follows:
$$
<\Phi(U)F^{[\gamma]}\,,\omega^{[\alpha]}(T)E_{[\alpha]}>
=<\Phi(U)\otimes F^{[\gamma]}\,,\Delta(\omega^{[\alpha]}(T))
(1\otimes E_{[\alpha]})>
$$
$$
=<\Phi(U)\,,\omega_{(1)}^{[\alpha]}(T)><F^{[\gamma]}\,,
\omega_{(2)}^{[\alpha]}(T)E_{[\alpha]}>
$$
$$
=<\Phi(U)\,,\omega_{(1)}^{[\alpha]}(T)><\Delta(F^{[\gamma]})\,,
E_{[\alpha]}\otimes\omega_{(2)}^{[\alpha]}(T)>
$$
$$
=<\Phi(U)\,,\omega_{(1)}^{[\alpha]}(T)>([N!]_B)_{[\alpha]}^{[\gamma]}
<1\,,\omega_{(2)}^{[\alpha]}(T)>
$$
\begin{equation}
\label{2.2}
=<\Phi(U)\,,([N!]_B)_{[\alpha]}^{[\gamma]}\omega^{[\alpha]}(T)>=0\,,
\end{equation}
where $\omega_{(1)}\otimes\omega_{(2)} := \Delta(\omega )$, and the
formula (\ref{1.17}) in the form
$ <F^{[\gamma]}\,,E_{[\alpha]}>=([N!]_B)_{[\alpha]}^{[\gamma]} $
has been used.
The equality (\ref{2.2}) must be fulfilled for arbitrary $\Phi(U)$.
Consequently,
due to non-degeneracy of the $T-U$-pairing,
the relation
\begin{equation}
\label{2.3}
[N!]_B\omega(T)=0
\end{equation}
is a necessary and sufficient condition for the relator (\ref{2.1})
to lie in the kernel of the bilinear form $<\,,\,>$.

\vspace{.1cm}

{\bf 10.}
We will finally show that the (remaining) null bi-ideals
are already generated by the elements (\ref{1.9}) and (\ref{1.9a})
satisfying (\ref{1.19}) and (\ref{1.19a}), respectively.
In other words,
considering the $T$-dependent case (\ref{2.1}),(\ref{2.3}) does not
enlarge the set of relators (\ref{1.9}) obtained with the use of all
solutions of the $T$-independent condition (\ref{1.19}).
Consider an element
$ \omega(T) \in V \otimes Y_-(T) $
satisfying (\ref{2.3}),
where $ V $ is the vector space associated to the multi-indices
$ [\alpha] $ appearing above
and $ Y_-(T) $ is the subalgebra in $ Y_- $ generated by the $ T $
generators.
Now the linear map
$ [N!]_B : V \otimes Y_-(T) \rightarrow V \otimes Y_-(T) $
is induced by
$ [N!]_B : V \rightarrow V, $
such that the kernel of the first map is
$ K \otimes Y_-(T), $
$ K $ being the kernel of the second one.
So, given all solutions of eq. (\ref{1.19}) and the corresponding
$q$-Serre relators (\ref{1.9}), we do not enlarge the null ideal
by considering
also the $T$-dependent case (\ref{2.3}).

\vspace{.1cm}

{\bf 11}.
Let us summarize the main steps in the construction of
bialgebras (and Hopf algebras) using the method proposed here.
We begin with
specifying any invertible matrix solution $R$ of the Yang-Baxter
equation
(for introducing an antipode, an additional `skew' invertibility
\cite{Vl1,FRT} of $R$ would be required).
Then we build bialgebras
(\ref{1}),(\ref{2})
(which may, under suitable conditions,
be extended to Hopf algebras \cite{Vl1})
being paired via (\ref{3}),(\ref{6}).
The double crossproduct construction is applied to these bialgebras in
order to obtain the `cross-multiplication' relations (\ref{ML1}).
Now we find `pure $T,U$' null ideals (they are
specific to a given $R$),
extend them as described in Sect. 6 and divide
out by the resulting bi-ideals.
The next step is to solve eq. (\ref{1.19}),
find all elements of type (\ref{1.9}) and their $F$-analogs (\ref{1.9a})
in
explicit form, and then equate them to zero.
This yields (if any solutions
of (\ref{1.19}) exist for a given $R$) a set of generalized $q$-Serre
relations which, for appropriate $R$, coincide with known $q$-Serre
relations of $U_q({\bf g})$.
Several examples of this procedure were presented in our
papers~\cite{Lud,Lud2,Vl2}.

\vspace{.1cm}

\noindent
{\small {\bf Acknowledgements}. This work was supported in part by RFBR
grant 97-01-01041. We thank S.\,Majid for many useful comments}.

\end{document}